\documentclass[11pt]{article}
\usepackage[colorlinks=true,
linkcolor=webgreen,
filecolor=webbrown,
citecolor=webgreen]{hyperref}

\usepackage{amssymb,psfig,epsfig}

\definecolor{webgreen}{rgb}{0,.5,0}
\definecolor{webbrown}{rgb}{.6,0,0}

\newtheorem{lemma}{Lemma}

\newcommand{\hsp}{\hspace*{\parindent}}
\newcommand{\eqn}[1]{(\ref{#1})}
\newcommand{\eeq}{\end{equation}}
\newcommand{\beql}[1]{\begin{equation}\label{#1}}
\newcommand{\bsq}{{\vrule height .9ex width .8ex depth -.1ex }}
\newcommand{\La}{\Lambda}
\newcommand{\af}{\alpha}
\newcommand{\sV}{{\mathcal V}}
\newcommand{\RR}{{\mathbb R}}
\newcommand{\HH}{{\mathbb H}}
\newcommand{\ZZ}{{\mathbb Z}}
\newcommand{\sP}{{\mathcal P}}
\newcommand{\sE}{{\mathcal E}}

\makeatletter
\def\@sect#1#2#3#4#5#6[#7]#8{\ifnum #2>\c@secnumdepth
     \def\@svsec{}\else
     \refstepcounter{#1}\edef\@svsec{\csname the#1\endcsname.\hskip .75em }\fi
     \@tempskipa #5\relax
      \ifdim \@tempskipa>\z@
        \begingroup #6\relax
          \@hangfrom{\hskip #3\relax\@svsec}{\interlinepenalty \@M #8\par}%
        \endgroup
       \csname #1mark\endcsname{#7}\addcontentsline
         {toc}{#1}{\ifnum #2>\c@secnumdepth \else
                      \protect\numberline{\csname the#1\endcsname}\fi
                    #7}\else
        \def\@svsechd{#6\hskip #3\@svsec #8\csname #1mark\endcsname
                      {#7}\addcontentsline
                           {toc}{#1}{\ifnum #2>\c@secnumdepth \else
                             \protect\numberline{\csname the#1\endcsname}\fi
                       #7}}\fi
     \@xsect{#5}}
\def\@begintheorem#1#2{\it \trivlist \item[\hskip \labelsep{\bf #1\ #2.}]}

\def\section{\@startsection {section}{1}{\z@}{-3.5ex plus -1ex minus 
 -.2ex}{2.3ex plus .2ex}{\normalsize\bf}}
\makeatother
\makeatletter
\def\subsection{\@startsection {subsection}{1}{\z@}{-3.5ex plus -1ex minus
 -.2ex}{2.3ex plus .2ex}{\normalsize\bf}}

\makeatother

\begin{document}
\begin{center}
{\large {\bf Quantizing Using Lattice Intersections}} \\
\vspace{1.5\baselineskip}
N. J. A. Sloane \\
Information Sciences Research Center \\
AT\&T Shannon Lab \\
Florham Park, New Jersey 07932-0971 \\
Email: \href{mailto:njas@research.att.com}{njas@research.att.com} \\ [+.25in]
B. Beferull-Lozano\footnotemark \\
Department of Electrical Engineering-Systems \\
University of Southern California \\
Los Angeles, CA 90089--2564 \\
Email: \href{mailto:beferull@sipi.usc.edu}{beferull@sipi.usc.edu} \\
\vspace{1.5\baselineskip}
August 27, 2001. Revised Feb 6, 2002. References updated July 17, 2002.\\
\vspace{1.5\baselineskip}
\begin{tabular}{l}
Dedicated to Eli Goodman and Ricky Pollack on \\
the occasion of their 66th/67th birthdays \\
\end{tabular} \\
\vspace{1.5\baselineskip}
{\bf ABSTRACT}
\vspace{.5\baselineskip}
\end{center}
\setlength{\baselineskip}{1.5\baselineskip}

\footnotetext[1]{This work was carried out during B. Beferull-Lozano's visit to AT\&T Labs in the summer of 2001.}

The usual quantizer based on an $n$-dimensional lattice $\Lambda$ maps a point $x \in \RR^n$ to a closest lattice point.
Suppose $\Lambda$ is the intersection of lattices $\La_1, \ldots, \La_r$.
Then one may instead combine the information obtained by simultaneously
quantizing $x$ with respect to each of the $\La_i$.
This corresponds to decomposing $\RR^n$ into a honeycomb of cells which are the intersections of the Voronoi cells for the $\La_i$, and identifying the cell to which $x$ belongs.
This paper shows how to write several standard lattices (the face-centered and body-centered cubic lattices, the root lattices $D_4$, $E_6^\ast$, $E_8$, the Coxeter-Todd, Barnes-Wall and Leech lattices, etc.) in a canonical way as intersections of a small number of simpler, decomposable, lattices.
The cells of the honeycombs are given explicitly and the mean squared quantizing error calculated in the cases when the intersection lattice is the face-centered or body-centered cubic lattice or the lattice $D_4$.

\section{Introduction}
\hsp
An $n$-dimensional lattice $\La$ determines a partition of $\RR^n$ into cells which are congruent copies of the Voronoi cell containing the origin.
A vector quantizer based on $\La$ (\cite{SPLAG}, \cite{GG92}, \cite{GrNe98}) is a map $Q_\La: \RR^n \to \RR^n$ which maps $x \in \RR^n$ to the closest lattice
point, i.e. to the lattice point at the center of the Voronoi cell containing $x$.
(In the case of a tie, one of the closest lattice points is chosen at random.)

Suppose $\La$ is the intersection of $n$-dimensional lattices $\La_1 , \ldots, \La_r$.
A recent paper \cite{BO1} discusses (among other things) the possibility of replacing the quantizer $Q_\La$ by a ``multiple description quantizer'' which simultaneously quantizes with respect to each of the $\La_i$, i.e. computes
\beql{EqQ1}
(Q_{\La_1} (x), Q_{\La_2} (x), \ldots, Q_{\La_r} (x)) ~.
\eeq

This gives rise to a different partition of $\RR^n$:
the cells are now the intersections of the Voronoi cells of the individual $\La_i$, and \eqn{EqQ1} specifies the cell to which $x$ belongs.

There are several reasons for investigating quantizers of this type.

(i) If the $\La_i$ are simpler lattices than $\La$ then \eqn{EqQ1} may be easier to compute than $Q_\La (x)$.

(ii) Perhaps this approach will lead to new insights on the so-far
intractable problem of finding good lattice quantizers in high dimensions (cf. \cite[Chap. 2]{SPLAG}).
Even in 24 dimensions the best lattice quantizer presently known, the Leech lattice, is quite complicated to analyze and to implement --- its Voronoi cell has 16969680 faces and over $10^{21}$ vertices
(\cite[Chaps. 21, 22, 23, 25]{SPLAG}, \cite{BeeSh92}, \cite{Vard95}, \cite{VaBe93}).

(iii) The individual $Q_{\La_i} (x)$ could be communicated over separate channels;
in the event of one or more channels failing a reasonably good approximation to $x$ will still be obtained.
Other multiple description quantizers have recently been studied in, for example, \cite{DSV}, \cite{VSS}.

(iv) Initially it seemed possible that this approach might lead to quantizers with a lower mean squared error that can be obtained from lattice quantizers.
The investigations reported here now make this unlikely (but see the remarks in the final section).

However, our focus here is not on applications but on the theoretical aspects of these quantizers, and in particular on two questions:
which lattices have a ``nice'' description as the intersection of simpler lattices, and what are the associated ``honeycombs,'' or
decompositions of space into cells?

\begin{figure}[htb]
\centerline{\epsfig{file=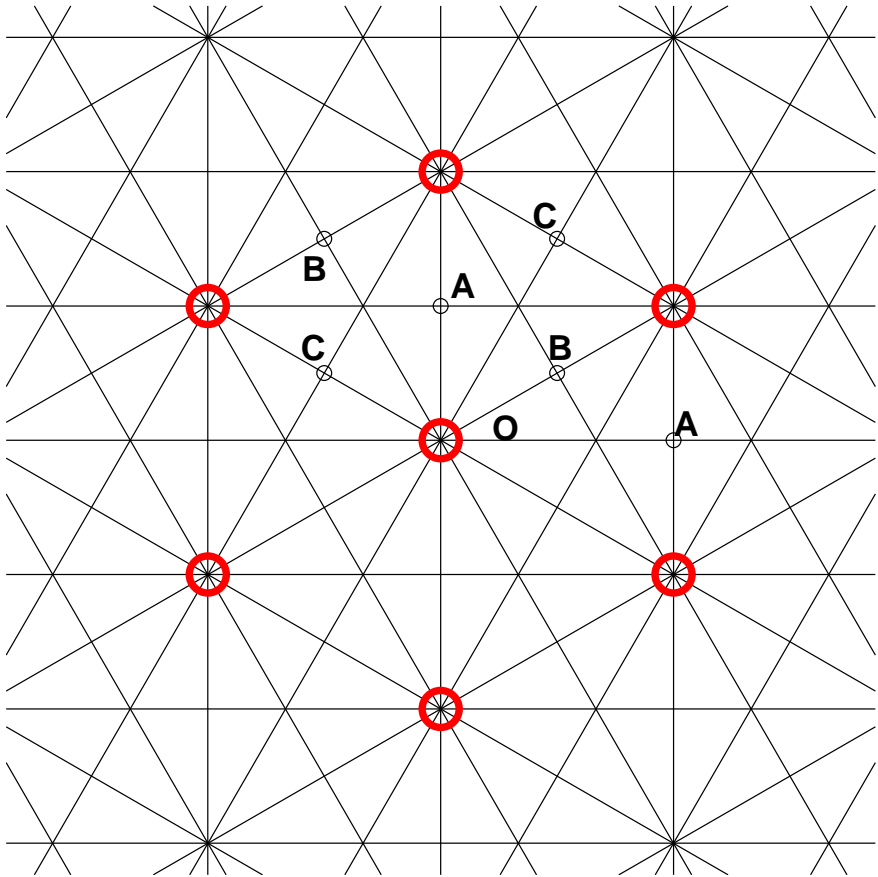,width=4.6in}}
\caption{The hexagonal lattice (heavy circles) as the intersection of three rectangular lattices (spanned by the vectors OA, OB and OC resp.).}
\label{A1}
\end{figure}

Reference \cite{BO1} gives one very appealing example, shown here in Figure \ref{A1}.
The familiar planar hexagonal lattice $A_2$ (large circles) can be obtained as the intersection of three rectangular lattices,
all rotations of each other by $120^\circ$:
these are the lattices generated respectively by the two vectors OA, the two vectors OB and the two vectors OC.
If we use generator matrices to specify these lattices (the rows span the lattices) then the three rectangular lattices $\La_1$, $\La_2$, $\La_3$ have generator matrices
\beql{EqA1}
\left[
\begin{array}{cc}
\sqrt{3} & 0 \\
0 & 1
\end{array}
\right] \,, \qquad
\left[\begin{array}{cc}
\frac{\sqrt{3}}{2} & \frac{1}{2} \\
\frac{- \sqrt{3}}{2} & \frac{3}{2}
\end{array}
\right] \,, \qquad
\left[
\begin{array}{cc}
\frac{- \sqrt{3}}{2} & \frac{1}{2} \\
\frac{\sqrt{3}}{2} & \frac{3}{2}
\end{array}
\right]
\eeq
and their intersection has generator matrix
\beql{EqA2}
\left[ \begin{array}{cc}
\sqrt{3} & 1 \\
0 & 2 \end{array}
\right] \,,
\eeq
which is indeed a copy of the $A_2$ lattice.

\begin{figure}[htb]
\centerline{\epsfig{file=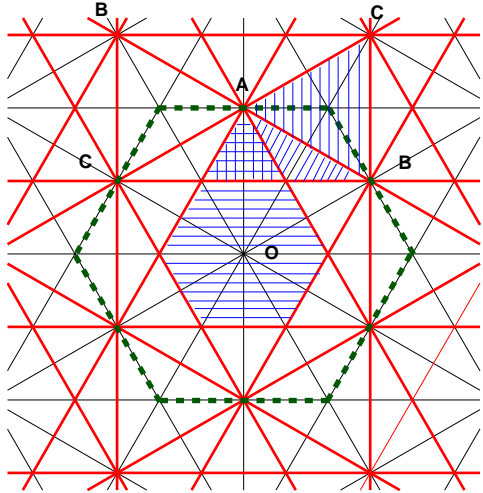,width=2.5in}}
\caption{Honeycomb associated with Voronoi cells in Fig. \ref{A1}.
[Note to copy editor: Figures \ref{A1} and \ref{A2} should be scaled so that the distance from 0 to the point labeled A in Fig. \ref{A2} is the same as that distance in Fig. \ref{A1}.]}
\label{A2}
\end{figure}

In this case the honeycomb contains four kinds of cells, as shown by the heavy (solid) lines in Fig. \ref{A2}.
For example,
if the point $x$ being quantized is close to $(0,0)$, then $Q_{\La_1} (x) = Q_{\La_2} (x) = Q_{\La_3} (x) =(0,0)$, and the cell containing $x$ is the intersection of the Voronoi cells containing $(0,0)$ of the three $\La_i$.
This is the horizontally shaded hexagon in Fig. \ref{A2}.
Just to the North of this hexagon the cell is the intersection of the Voronoi cell for $\La_1$ that contains $(0,1)$ with the Voronoi cells at $(0,0)$ for $\La_2$ and $\La_3$:
this is the small cross-hatched equilateral triangle.
There are two further cells that are obtained in a similar manner:
the diagonally shaded isosceles triangle and the larger vertically shaded equilateral
triangle.
This example will be discussed further in Section 4.

The present paper was prompted by the question: how can Figures \ref{A1} and \ref{A2} be generalized to higher dimensions?

Section 2 gives some terminology and defines the normalized mean squared error $G$ that will be used when comparing quantizers.
Section 3 describes some general constructions for writing a lattice as an intersection and gives a number of examples,
including the body-centered cubic (bcc) and face-centered cubic (fcc) lattices
$D_3^\ast$ and $D_3$,
the root lattices $D_4$, $E_6^\ast$, $E_8$, the Coxeter-Todd lattice $K_{12}$,
the Barnes-Wall lattices $BW_n$ and the Leech lattice $\La_{24}$.
We focused attention on these lattices because
$A_2$, $D_3^\ast$, $D_4$, $E_6^\ast$, $E_8$, $K_{12}$, $BW_{16}$ and $\La_{24}$
are the best quantizers currently known\footnote{Assuming always that the variable to be
quantized is uniformly distributed over a large ball in $\RR^n$.}
in their dimensions \cite{SPLAG}.
In fact $A_2$ is optimal among all two-dimensional quantizers \cite{Fej8},
and $D_3^\ast$ is optimal among three-dimensional lattice quantizers \cite{BS1}.

In Section 4 we reexamine the quantizer of Figures \ref{A1} and \ref{A2}.
Then Sections 5, 6 and 7 determine the honeycombs and mean squared errors for the bcc, fcc and $D_4$ lattices, respectively.
The final section contains some conclusions and comments.

Table~\ref{TS} summarizes the main decompositions mentioned in this paper.

\begin{table}[htb]
$$
\begin{array}{cccccccc}
\mbox{Lattice} & \mbox{Copies} & \mbox{Component} & \mbox{Sections} \\
& & \mbox{lattice} & \\ [+.05in]
A_2 & 3 & \mbox{rectangular} & 1, 3, 4 \\
A_3^\ast ~\mbox{(bcc)} & 3 & \mbox{``rectangular''} & 3, 5 \\
A_3 ~\mbox{(fcc)} & 4 & \mbox{``prismatic''} & 3, 6 \\
D_4 & 3 & \ZZ^4 & 3, 7 \\
E_6^\ast & 4 & A_2^\ast & 3 \\
E_8 & 15 & \ZZ^8 & 3 \\
E_8 & 10 & A_2^4 & 3 \\
E_8 & 5 & D_4^2 & 3 \\
K_{12} & 21 & A_2^6 & 3 \\
\mbox{Leech} & 4095 & \ZZ^{24} & 3 \\
BW_n & \prod_{j=1}^{m-1} (2^j -1 ) & \ZZ^n & 3 
\end{array}
$$
\caption{Summary of decompositions described in this paper.}
\label{TS}
\end{table}

We repeat that we do not propose that the quantizers described in this paper will have practical applications, 
nor that they are simpler than the usual lattice quantizers based on the intersection lattices.
Our interest is in the geometry of the honeycombs
associated with the quantizers.
These do not seem to have been studied before\footnote{Of course
there is an extensive literature dealing with the Voronoi and Delaunay
cell decompositions associated with various lattices --- see for example \cite{CSV}, \cite{CSF}, \cite{CSC}, \cite{SPLAG}, \cite{Cox}.}
--- they are not mentioned for example in
Okabe, Boots and Sugihara \cite{OBS}, nor in Wells \cite{Wells},
even in the chapter on {\em Interpenetrating three-dimensional nets}.
There has been considerable interest among physicists in recent years in soap froths, both in connection with the study of colloids, dendritic polymers, etc., and in the Kelvin problem of finding the minimal-area honeycomb in $\RR^3$
(\cite{Riv}, \cite{WP94}, \cite{ZK1}).
So it may be of interest to see a class of honeycombs that have a purely geometrical construction.

We end this section by describing how the cells of the honeycombs in Sections 5--7 were found.
Many cells could be found by elementary geometrical
reasoning. But in complicated cases we carried out some or all of the following steps.

(i) To find the cell containing a point $x \in \RR^3$ or $\RR^4$ we
first quantized $x$ using each of the lattices $\La_i$ in turn.
For each $\La_i$, we determined the Voronoi cell $C_i$ containing $x$, or, more precisely, the equations to the hyperplanes bounding $C_i$.

(ii) Linear programming (in MATLAB) was then
used to determine the vertices of the cell containing $x$.
We let $w$ range over a set of 130 points on a sphere centered at $x$ (taken
from the tables of spherical codes in \cite{Har2a}), and, for each $w$, we maximized
the inner product $(w, z)$ subject to the constraints that $z$ lie in the polytope formed
by the intersection of {\em all}
the hyperplanes found in (i).
Any such solution $z$ is a vertex of the cell, and
since the $w$'s are essentially random,
the $130$ solutions should include all the vertices of the cell.

(iii) The convex hull program QHULL \cite{QHULL2}, \cite{QHULL} was used
to find the convex hull of these vertices.

At this point we have candidates for all the cells in the honeycomb.
Since it is theoretically possible (although unlikely)
that step (i) might have failed to find all the vertices
of a cell, we now verified by hand that the cells fitted together
to form a proper tiling of the space.

(iv) The XGobi program \cite{XGobi} for displaying
multi-dimensional data was used to help visualize the cells
and their neighbors.

To compute the volumes and second moments of the cells we decomposed the cells
into simplices and used the formulae in \cite{CSV} and \cite[Chap. 21]{SPLAG}.

\section{Notation}
\hsp
We assume the $n$-dimensional lattice $\La$ is the intersection of $n$-dimensional lattices $\La_1, \ldots, \La_r$.
The intersections of the Voronoi cells for the $\La_i$ partition $\RR^n$ into a tesselation or {\em honeycomb}.
Let $\sP_1, \ldots, \sP_k$ be representatives for the different polytopes or {\em cells} that appear in the honeycomb.

Let $p_i$ $(i=1, \ldots, k)$ be the probability that a 
randomly chosen point in $\RR^n$ (uniformly distributed over a very large ball, say)
belongs to $\sP_i$, and let $N_i = p_i V / V_i$,
where $V_i$ is the volume of $\sP_i$ and $V= \sqrt{\det \La}$ is the volume of a fundamental region or Voronoi cell for $\La$.
Let $\sV$ be the particular Voronoi cell for $\La$ that contains the origin.
Then the honeycomb is periodic with ``tile'' equal to the part lying in $\sV$, and there are $N_1$ cells of type $\sP_1$ per copy of $\sV$, $N_2$ cells of type $\sP_2$, etc.
Also
\beql{EqU0}
V= N_1 V_1 + \cdots + N_k V_k ~.
\eeq

We assume that a point that falls into a cell of type $\sP_i$ is quantized as the centroid $c_i$ of that cell,
in which case the mean squared error is
$$U_i = \int_{\sP_i} \| x-c_i \|^2 dx ~.$$
The mean squared error for the full quantizer is then
\beql{EqU1}
U = \sum_{i=1}^k N_i U_i = V \sum_{i=1}^k \frac{p_i U_i}{V_i} ~.
\eeq

We now derive the expression that we will use as a measure of the normalized mean squared error of this quantizer.
Suppose we are quantizing a random variable $X \in \RR^n$, with differential
entropy per dimension $h(X)$,
and whose support contains a large number of points of $\La$.

A general theorem of Zador
about vector quantizers
(\cite{Z1}, \cite{Z2}, \cite[Eq. (20)]{GrNe98})
implies that in this ``high-rate'' case
the average mean squared error per dimension $U/(nV)$ can be approximated by
\beql{EqU2}
\frac{U}{nV} \approx G \, 2^{2(h(X) -R)} ~,
\eeq
where $R$ bits/symbol is the quantizing rate and $G$ depends on the positions of the quantizing points but is independent of $X$.

Since $G$ does not depend on the distribution of $X$,
we may choose any convenient distribution in order to calculate $G$,
and we assume that $X$ is uniformly distributed over a large region of $\RR^n$,
or, equivalently, that $X$ is uniformly distributed over $V$.
Then
\beql{EqU3}
h(X) = \frac{1}{n} \log_2 V, ~~~ {\rm so} ~~~
2^{2h(X)} = V^{2/n} ~.
\eeq

To calculate $R$, observe that we need $h(p_1, \ldots, p_k ) = - \sum_{i=1}^k p_i \log_2 p_i$
bits to specify the type of cell to which the quantized
point belongs, and a further
$\sum_{i=1}^k p_i \log_2 N_i = \sum_{i=1}^k p_i \log_2 (p_i V/V_i )$ bits
to specify the particular one of the $N_i$ cells of that type.
This requires a total of $\log_2 V - \sum_{i=1}^k p_i \log_2 V_i$ bits, and
then $R$ is this quantity divided by $n$, so that
\beql{EqU4}
2^{-2R} = V^{-2/n} \prod_{k=1}^k V_i^{2p_i /n} ~.
\eeq
From \eqn{EqU1}--\eqn{EqU4} we obtain
\beql{EqU5}
G = \frac{\sum_{i=1}^k \frac{p_i U_i}{V_i}}{n \left( \prod_{i=1}^k V_i^{p_i}\right)^{\frac{2}{n}}} ~,
\eeq
the normalized mean squared error per dimension, which we take as our figure of merit
for a quantizer.
The numerator of \eqn{EqU5} is equal to $U/V$ (see \eqn{EqU1}).
Note that when there is only one kind of cell \eqn{EqU5} reduces to the familiar formula
\beql{EqU6}
G = \frac{U}{nV^{1+ \frac{2}{n}}}
\eeq
for a lattice quantizer (\cite{CSV}, \cite{SPLAG}, \cite{GG92}).

Incidentally, a different expression from \eqn{EqU5} for the figure of merit was used in a recent paper of Kashyap and Neuhoff \cite{KaNe}.
By defining the rate of the quantizer in a different way they end up with a different expression in the denominator.
However, we believe our formula gives a fairer comparison.
This point is discussed in more detail in \cite{SV}.

\paragraph{Further terminology.}
Let $\La$ be a lattice in $\RR^n$.
The dual lattice will be denoted by $\La^\ast$.
The {\em norm} of a vector $x \in \RR^n$ is its squared length $(x, x)$.
A {\em similarity} $\sigma$ is a linear map from $\RR^n$ to $\RR^n$ such that there is a real number $N$ with $(\sigma x, \sigma y) = N (x, y)$ for $x,y \in \RR^n$.
If $\La$ and $M$ are similar lattices we write $\La \cong M$.
A lattice $\La$ is  said to be {\em $N$-modular}
if $\La \cong \La^\ast$ under a similarity that multiplies norms by $N$
(cf. \cite{Queb}).
For example, the root lattices $\ZZ^n$ and $E_8$ are 1-modular,
$A_2$ is 3-modular, and $D_4$ is 2-modular.
We write $\langle ~ \La_1, \La_2, \ldots ~  \rangle$ for the lattice
generated by the union of the lattices $\La_1, \La_2 , \ldots$.
Two lattices or polytopes are {\em congruent} if one can be mapped to the other by an element of $SO(n)$.

\section{Writing a lattice as an intersection}
\hsp
Let $\La$ be a lattice in $\RR^n$.
We wish to write $\La$ as an intersection
\beql{EqD0}
\La = \La_1 \cap \La_2 \cap \cdots \cap \La_r
\eeq
where $r$ is small
and the $\La_i$ are pairwise congruent and as ``simple'' as possible.
Ideally we would like each $\La_i$ to be a direct sum of congruent copies of a fixed low-dimensional lattice $K$ such as $\ZZ$, $A_2$ or $D_4$, but this is not always possible.
In the example shown
in Fig.~\ref{A1}, for instance, the $\La_i$ are rectangular rather than square lattices.
If this is not possible, we ask that the $\La_i$ be decomposable into a direct sum of as many congruent low-dimensional sublattices as possible.

If we were going to investigate the honeycombs associated with higher-dimensional intersections such as those for $E_6^\ast$ or $E_8$, we would impose an additional formal requirement that the $\La_i$ form an orbit under some subgroup of the automorphism group $Aut (\La )$, in order to guarantee that the honeycomb be symmetric.
However, for the low-dimensional examples studied in Sections 5--7, we were able to achieve this symmetry by using the natural decompositions, without introducing the machinery of group theory.

This section describes some general methods for finding intersections.
We will make use of the following version of de~Morgan's law:

\begin{lemma}\label{DM}
If $\La_1, \ldots, \La_r$ are lattices in $\RR^n$ then
$$\La_1 \cap \cdots \cap \La_r ~=~ \langle ~ \La_1^\ast , \ldots, \La_r^\ast ~ \rangle^\ast ~.
$$
\end{lemma}

\paragraph{Proof.}
See for example \cite[Section 82F]{OM}.~~~$\bsq$

\paragraph{Method 1: Partitioning the minimal vectors for $\La^\ast$.}
For the first method, suppose the dual lattice $\La^\ast$ is generated by its minimal vectors (i.e. the vectors of minimal nonzero norm), and let $K$ be one of $\ZZ$, $A_2$ or $D_4$.
Let $K$ have dimension $\kappa$ and be $N$-modular.
Then we attempt to partition the minimal vectors of $\La^\ast$ into a number of copies of the minimal vectors of a suitably rescaled version of $K^{n/\kappa}$.

If this can be done, let $\La_i^\ast$ be the lattice generated by the $i$th part of the partition.
The hypotheses guarantee that $\La_i^\ast \cong \La_i \cong K^{n/\kappa}$,
and, by the Lemma,
$$\La = \La_1 \cap \cdots \cap \La_r ~,$$
where $r =$ number of minimal vectors of $\La^\ast$ divided by number of minimal vectors of $K$.

Whether this partitioning is possible is an interesting question in its own right.
For example, can the 240 minimal vectors of $E_8$ be partitioned into 15 copies of the minimal vectors of (a scaled version of) $\ZZ^8$, i.e. into 15 coordinate frames\footnote{A {\em coordinate frame} in $\RR^n$ is a set of $2n$ vectors $\pm v_1 , \ldots, \pm v_n$ with $(v_i, v_i) =$ a constant, $(v_i, v_j) =0$ if $i \neq j$.}
or into 10 copies of the minimal vectors of $A_2^4$?
Partial answers are given below.

We now give a number of examples, beginning with the case when $K = \ZZ$.

\paragraph{${\mathbf D_4}$.}
The lattice $D_4$ may be taken to have generator matrix
\beql{EqD1}
\left[\begin{array}{cccc}
2 & 0 & 0 & 0 \\
0 & 2 & 0 & 0 \\
0 & 0 & 2 & 0 \\
1 & 1 & 1 & 1
\end{array}
\right]
\eeq
and then the 24 minimal vectors consist of eight of the form $(\pm 2,0,0,0)$ and
16 of the form $(\pm 1, \pm 1 , \pm 1, \pm 1 )$.
$D_4$ is 2-modular and is generated by its minimal vectors.
The minimal vectors may be partitioned into three coordinate frames, consisting of $\pm 1$ times the rows of each of the matrices
\beql{EqD2}
\left[\begin{array}{cccc}
2 & 0 & 0 & 0 \\
0 & 2 & 0 & 0 \\
0 & 0 & 2 & 0 \\
0 & 0 & 0 & 2
\end{array}\right] \,, \quad
\left[\begin{array}{cccc}
+1 & +1 & +1 & +1 \\
+1 & -1 & +1 & -1 \\
+1 & -1 & -1 & +1 \\
+1 & +1 & -1 & -1
\end{array}
\right] \,, \quad
\left[\begin{array}{cccc}
-1 & +1 & +1 & +1 \\
-1 & -1 & +1 & -1 \\
-1 & -1 & -1 & +1 \\
-1 & +1 & -1 & -1
\end{array}
\right] \,.
\eeq
After applying the Lemma and rescaling, we conclude that if $\La_1$, $\La_2$, $\La_3$ $(\cong \ZZ^4 )$ have the generator matrices given in \eqn{EqD2}, then
\beql{EqD3}
\La_1 \cap \La_2 \cap \La_3 = \La \,,
\eeq
where $\La$ has generator matrix
\beql{EqD4}
\left[
\begin{array}{cccc}
4 & 0 & 0 & 0 \\
2 & 2 & 0 & 0 \\
2 & 0 & 2 & 0 \\
2 & 0 & 0 & 2
\end{array}
\right] \,,
\eeq
and is another version of $D_4$ on the scale at which its minimal norm is 8.
Equation (\ref{EqD3}) may also be verified directly, without appealing to the Lemma.
The group generated by the second matrix in \eqn{EqD2} and
${\rm diag} \{-1, +1, +1, +1\}$ is a symmetric group $S_3$ permuting the $\La_i$;
it is also a subgroup of $Aut (\La )$.
The honeycomb for this example is studied in Section 7.

\paragraph{Barnes-Wall lattices.}
The preceding example can be generalized using orthogonal spreads.
Let $BW_n$ $(n=2^m$, $m=1,2, \ldots )$ denote the $n$-dimensional Barnes-Wall lattice (\cite{SPLAG}, \cite{cliff1}, \cite{BW}).
In particular, $BW_2 \cong \ZZ^2$, $BW_4 \cong D_4$, $BW_8 \cong E_8$.
For $n \neq 8$, $BW_n$ is 2-modular, while as already mentioned $BW_8$ is 1-modular.

It is known that the minimal vectors of $BW_n$ may be partitioned into
$\prod_{j=1}^{m-1} (2^j -1 )$ coordinate frames, which are transitively
permuted by symmetries of $Aut (BW_n )$.
This is a consequence of the existence of an orthogonal spread in the orthogonal vector space $\Omega^+ (2m,2)$ of maximal Witt index
(\cite{CCKS}, \cite{grass3}, \cite{cliff1}).

It follows that $BW_n$ can be written as the intersection of $\prod_{j=1}^{m-1} (2^j -1 )$ copies of $\ZZ^n$.
In particular, $E_8$ is the intersection of 15 copies
of $\ZZ^8$.
An explicit method for constructing such an intersection for $E_8$ is given below.

Intersections of smaller numbers of lattices are possible, although they are less symmetric and therefore less satisfactory.
For example in \eqn{EqD3} it is also true that $\La_1 \cap \La_2 = \La \cong D_4$.
Similarly, $E_8$ is (up to a similarity) the intersection of $2\ZZ^8$ and the lattice (similar to $2\ZZ^8$) with generator matrix
$$
\left[
\begin{array}{rrrrrrrr}
1 & 1 & 1 & 1 & 0 & 0 & 0 & 0 \\
0 & 0 & 0 & 0 & 1 & 1 & 1 & 1 \\
1 & - 1 & 0 & 0 & 1 & -1 & 0 & 0 \\
0 & 0 & 1 & -1 & 0 & 0 & 1 & -1 \\
1 & 0 & -1 & 0 & -1 & 0 & 1 & 0 \\
0 & 1 & 0 & -1 & 0 & -1 & 0 & 1 \\
1 & 0 & 0 & -1 & 0 & 1 & -1 & 0 \\
0 & 1 & -1 & 0 & 1 & 0 & 0 & -1
\end{array}
\right]
$$
But this representation of $E_8$ is in no way canonical, and the resulting
honeycomb is not interesting.

\paragraph{Eisenstein and Hurwitzian lattices.}
Smaller intersections which {\em are} canonical can be obtained if we change $K$ from $\ZZ$ to $A_2$ or $D_4$.
For example, the minimal vectors of $E_8$ can be partitioned into 10 copies of the minimal vectors of $A_2^4$.
As in \cite{SPLAG}, let $\sE = \{a+b\omega : a,b \in \ZZ \}$, $\omega = e^{2 \pi i/3}$, denote the ring of Eisenstein integers.
The six units in $\sE$ are $\pm 1$, $\pm \omega$, $\pm \bar{\omega}$.
When regarded as a two-dimensional real lattice $\sE$ is similar to $A_2$.
As
an $\sE$-module, $E_8$ has generator matrix
\beql{EqE1}
\left[\begin{array}{rrrr}
\theta & 0 & 0 & 0 \\
0 & \theta & 0 & 0 \\
1 & 1 & 1 & 0 \\
0 & 1 & -1 & 1
\end{array}
\right] \,, \quad
\mbox{where} \quad \theta = \omega - \bar{\omega} \,.
\eeq
Inner products are computed using the Hermitian inner product $(u,v) = \sum u_i \bar{v}_i$.
See \cite[Chapters 2 and 7]{SPLAG} for further details.
The minimal vectors consist of 24 of the form $(u \theta , 0,0,0)$, where $u$ is a unit in $\sE$, and $8 \times 3^3 = 216$ which are congruent $\bmod~\theta$ to one of the eight nonzero codewords of the tetracode \cite[Chap. 3]{SPLAG}.

A partition of these 240 vectors into 10 copies of the minimal vectors of $A_2^4$ was found by graph coloring.
A graph was constructed with the 40 projectively distinct vectors as nodes and with edges corresponding to pairs of non-orthogonal vectors.
A coloring with 10 colors was then found with the help of a program supplied by David Johnson.
The ten copies of $\sE^4 \cong A_2^4$ are shown in Table \ref{TEA}.
\begin{table}[htb]
$$
\begin{array}{rrrr|rrrr|rrrr|rrrr}
\theta & 0 & 0 & 0 & 0 & \theta & 0 & 0 & 0 & 0 & \theta & 0 & 0 & 0 &0 & \theta \\
0 & 1 & -1 & 1 & 1 & 0 & -1 & -1 & 1 & -1 & 0 & 1 & 1 & 1 & 1 & 0 \\
0 & 1 & -1 & \omega & 1 & 0 & -\omega & -\bar{\omega} & 1 & -\omega & 0 & \bar{\omega} & 1 & \omega & \omega & 0 \\
0 & 1 & -\omega & \omega & 1 & 0 & -\omega & -\omega & 1 & -1 & 0 & \omega & 1 & 1 & \omega & 0 \\
0 & 1 & -\omega & \bar{\omega} & 1 & 0 & - \bar{\omega} & -1 & 1 & -\omega & 0 & 1 & 1 & \omega & \bar{\omega} & 0 \\
0 & 1 & -\omega & 1 & 1 & 0 & -1 & -\bar{\omega} & 1 & -\bar{\omega} & 0 & \bar{\omega} & 1 & \bar{\omega} & 1 & 0 
\end{array}
$$

\caption{Decomposition of minimal vectors of $E_8$ into ten copies of ${\mathcal E}^4 \stackrel{\sim}{=} A_2^4$.
Each row generates a copy of ${\mathcal E}^4 \stackrel{\sim}{=} A_2^4$.  The complex conjugates of the last four rows have been omitted.}
\label{TEA}
\end{table}
Only one from each complex conjugate pair is shown.
By applying the Lemma we obtain a representation of $E_8$ as an intersection of 10 copies of $A_2^4$.
In fact (since ${\mathcal E}$ is itself $3$-modular)
we may omit the final step of taking the duals of the lattices in Table \ref{TEA}.
Let $\La_1, \ldots, \La_{10}$ be the ten versions of $A_2^4$ generated by the rows of Table \ref{TEA} and their complex conjugates.
Then their intersection is easily seen to be the version of $E_8$ with generator matrix $\theta$ times \eqn{EqE1}.
This decomposition is probably not unique, and it would be nice to know which version has the largest symmetry group.

Again just two lattices suffice:
$E_8$ is also the intersection of the first two lattices in Table \ref{TEA}.

We may also write $E_8$ as the intersection of five copies of $D_4^2$.
For this we regard $E_8$ as a 2-dimensional module over the ring $\HH \cong D_4$ of Hurwitzian quaternions \cite[p. 55]{SPLAG}.
The five copies of $D_4^2$ have generator matrices
\beql{EqE2}
\left[
\begin{array}{cc}
1+i & 0 \\
0 & 1+i
\end{array}\right] ,~~
\left[\begin{array}{cr}
1 & 1 \\ 1 & -1
\end{array}\right] ,~~
\left[\begin{array}{cr}
1 & i \\ 1 & -i
\end{array}\right] ,~~
\left[\begin{array}{cr}
1 & j \\ 1 & -j
\end{array}\right], ~~
\left[\begin{array}{cr}
1 & k \\ 1 & -k
\end{array}\right] \,.
\eeq

$E_6^\ast$ may be written as the intersection of four copies of $A_2^3$, with generator matrices
\beql{EqE3}
\left[\begin{array}{ccc}
\theta & 0 & 0 \\
0 & \theta & 0 \\
0 & 0 & \theta
\end{array}\right] ,~~
\left[\begin{array}{ccc}
1 & 1 & 1 \\
1 & \omega & \bar{\omega} \\
1 & \bar{\omega} & \omega
\end{array}\right],~~
\left[\begin{array}{ccc}
1 & 1 & \omega \\
1 & \omega & 1 \\
1 & \bar{\omega} & \bar{\omega}
\end{array}\right] , ~~
\left[\begin{array}{ccc}
1 & 1 & \bar{\omega} \\
1 & \omega & \omega \\
1 & \bar{\omega} & 1
\end{array}\right] \,.
\eeq
This was found by partitioning the 72 minimal vectors of $E_6$ (by hand) into four copies of the minimal vectors of $A_2^3$ and using the Lemma.

\paragraph{Method 2: Congruence bases and norm-doubling maps.}
The second method is based on the observation that several well-known lattices $\La$ have the property that for some prime $\pi$, the vectors in some of the classes of $\La / \pi \La$ can be partitioned into coordinate frames.
For example, Conway's proof of the uniqueness of the Leech
lattice $\La_{24}$ \cite[Chap. 12]{SPLAG} considers the classes of $\La_{24} / 2\La_{24}$.
A consequence of the numerical identity
\begin{eqnarray}
\label{EqE4}
\frac{n_0}{1} +
\frac{n_4}{2} + \frac{n_6}{2} + \frac{n_8}{48}
& = &
\frac{1}{1} +
\frac{196560}{2} + \frac{16773120}{2} + \frac{398034000}{48} \nonumber \\
& = & 16777216 = 2^{24} \,,
\end{eqnarray}
where $n_j$ is the number of vectors in $\La_{24}$ of norm $j$, is that, for the classes of $\La_{24} / 2 \La_{24}$ in which the minimal norm is 8, the minimal vectors in the class form a coordinate frame or {\em congruence base}.
A similar property holds for the $D_4$, $E_8$, $K_{12}$ and other lattices
(cf. \cite{K12}).
This gives a representation of $\La_{24}$ as an intersection of
$398034000/48 = 8292375$ copies of $\ZZ^{24}$.
However, the following argument, due to J.~H. Conway (personal communication), shows that if the lattice has a suitable norm-doubling map (cf. \cite[p. 239]{SPLAG}, \cite{K12}) then we can also
partition the minimal vectors into coordinate frames and obtain a smaller intersection.

Suppose a lattice $\La \subseteq \RR^n$ has the structure of a free module over a ring $J$ with inner product $(~,~)$ (cf. \cite[Chap. 2]{SPLAG}).
In the present application $J$ will be either $\ZZ$ or $\sE$.
Let $a = n/ {\rm dim}_\RR J$.
Consider the classes of $\La / 2 \La$.
Suppose there is an integer $m$ with the property that each class either contains no vectors of norm $2m$, or else all the vectors of norm $2m$ in the class can be partitioned into sets of $2a$ vectors
$\pm v_1$, $\pm v_2 , \ldots, \pm v_a$ where $(v_i, v_j ) =0$ if $i \neq j$.

Suppose in addition there is a {\em norm-doubling} map $T$, a similarity from $\La$ into $\La$ such that $(Tu, Tu) = 2(u,u)$ for $u \in \La$, with the extra property that $2\La \subset T\La$.
Then we may conclude that the vectors of norm $m$ in $\La$ may also be partitioned into sets of $2a$ mutually orthogonal vectors.

To see this, let $u \in \La$ have norm $m$.
Then $v = Tu$ has norm $2m$, and by the hypotheses is part of a coordinate frame $\pm v_1 , \ldots \pm v_a$, where
$v_i = v+ 2w_i$, say, with $w_1 =0$.
We can write $2w_i = Tw'_i$ for some $w'_i$, so $v_i = T(u+w'_i )$.
Since $T$ is a similarity, the set $\pm (u+ w'_i )$ is a coordinate frame containing $u$.

\subsection*{Examples}
\hsp
(i) $\La = D_4$ or $E_8$, $J = \ZZ$, $m=2$, $T=$ direct sum of
respectively 2 or 4 copies of $\left[ \begin{array}{cr} 1&1 \\ 1& -1 \end{array}\right]$, with $T^2 = 2I$.
The classes of $D_4/2D_4$ and $E_8 /2E_8$ and the associated congruence bases are given in \cite[Chap. 6]{SPLAG}.
The analogue of \eqn{EqE4} for $E_8$ reads
$$1 + \frac{240}{2}
+ \frac{2160}{16} =2^8 ~.
$$
We obtain decompositions of the minimal vectors of $D_4$ into three coordinate frames, as already seen in \eqn{EqD2}, and of the minimal vectors of $E_8$ into 15 coordinate frames as also discussed above.
To get an explicit decomposition in the latter case, note that a coordinate frame of norm 4 vectors has the form
$Tu + 2w_i = T(u + Tw_i )$.
So the coordinate frame of norm 2 vectors consists of the vectors of minimal norm in the translate $u+T\La$.
For $E_8$ these consist of seven sets of the form shown on the left in \eqn{EqE5} and eight of the form shown on the right:
\beql{EqE5}
\left[\begin{array}{cccccccc}
+ & 0 & + & 0 & 0 & 0 & 0 & 0 \\
+ & 0 & - & 0 & 0 & 0 & 0 & 0 \\
0 & + & 0 & + & 0 & 0 & 0 & 0 \\
0 & + & 0 & - & 0 & 0 & 0 & 0 \\
0 & 0 & 0 & 0 & + & 0 & + & 0 \\
0 & 0 & 0 & 0 & + & 0 & - & 0 \\
0 & 0 & 0 & 0 & 0 & + & 0 & + \\
0 & 0 & 0 & 0 & 0 & + & 0 & -
\end{array}
\right] \,, \qquad\frac{1}{2} \left[
\begin{array}{cccccccc}
+ & + & + & + & + & + & + & + \\
+ & + & + & + & - & - & - & - \\
+ & + & - & - & + & + & - & - \\
+ & + & - & - & - & - & + & + \\
+ & - & + & - & + & - & + & - \\
+ & - & + & - & - & + & - & + \\
+ & - & - & + & + & - & - & + \\
+ & - & - & + & - & + & + & -
\end{array}\right] \,.
\eeq
Then $E_8$ is also the intersection of the 15 copies of $\ZZ^8$ having these generator matrices.

(ii) $\La =$ Leech lattice $\La_{24}$, $J= \ZZ$, $m=4$, $T= (I+i)$, where $i \in Aut (\La_{24})$ is given 
in \cite[Fig. 6.7]{SPLAG}, and satisfies $i^2 =-I$, with $(I+i) (I-i) = T(I-i) =2I$.
The coordinates frames of norm 4 vectors consist of the 48 vectors of minimal norm in the translates $u+ (I-i) \La$ where $(u,u)=4$.
We obtain a decomposition of the minimal vectors of $\La_{24}$ into 4095 coordinate frames, and a representation of $\La_{24}$ as the intersection of 4095 copies of $\ZZ^{24}$.

We do not know if it is possible to write the Leech lattice as the intersection of 2730 copies of $A_2^{12}$.
Since 196560/96 is not an integer, there is no analogous decomposition as an intersection of copies of $D_4^6$.

(iii) $\La =$ Coxeter-Todd lattice $K_{12}$, a 6-dimensional $\sE$-module, $J= \sE$, $m=6$, $T=$ the map given in \cite[Eq. (43)]{K12}, with $T^2+ T+2 =0$.
The classes of $K_{12} /2K_{12}$ are given in \cite{K12}, and the analogue of \eqn{EqE4} reads
$$1 + \frac{756}{2} + \frac{4032}{1} + \frac{20412}{12} = 2^{12} ~.$$
The coordinate frames of norm 6 vectors consist of the 12 minimal vectors in the translates
\linebreak
$u+ (T+I)K_{12}$, $(u,u) =6$.
By combining these coordinate frames in sets of three, by taking the union of
the sets $\af \{ \pm v_1, \ldots, \pm v_{12} \}$ with $\af = 1, \omega$ and $\bar{\omega}$,
we obtain a decomposition of the minimal vectors of $K_{12}$ into 21 copies of the minimal vectors of $A_2^6$, and, via the Lemma,
a representation of $K_{12}$ as the intersection of 21 copies of $A_2^6$.
An explicit decomposition, not shown here, was found by the graph coloring method mentioned earlier.

\paragraph{Method 3: First principles.}
If the above methods fail, as they do for the bcc and fcc lattices,
we can always fall back on a direct attack from first principles.
The following method handles the hexagonal, bcc and fcc lattices in a unified manner.
We list the vectors of small norms in the lattice, and look for a partition of some subset
of these vectors which produces a small number of congruent, decomposable lattices whose
intersection is similar to the original lattice.

For the hexagonal lattice, which we take to be generated by
$(0,1)$ and $( - \frac{\sqrt{3}}{2} , \frac{1}{2} )$, there are six vectors
of norm 1, namely $(0, \pm 1)$, $( \pm \frac{\sqrt{3}}{2} , \pm \frac{1}{2} )$, and six of norm 3,
namely $( \pm \sqrt{3}, 0)$, $(\pm \frac{\sqrt{3}}{2} , \pm \frac{3}{2} )$.
Then \eqn{EqA1} is obtained by partitioning these 12 vectors into three sets of size 4.

For the bcc lattice generated by $(1,0,0)$, $(0,1,0)$,
$(\frac{1}{2}, \frac{1}{2}, \frac{1}{2} )$, the vectors of small norms are the following:
$$
\begin{array}{ccccccc}
\multicolumn{3}{c}{\mbox{shape}} & ~~~~~ & \mbox{norm} & ~~~~~ & \mbox{number} \\ [+.05in]
0 & 0 & 0 && 0 && 1 \\
\frac{1}{2} & \frac{1}{2} & \frac{1}{2} && \frac{3}{4} && 8 \\
1 & 0 & 0 && 1 && 6 \\
1 & 1 & 0 && 2 && 12
\end{array}
$$
We take the 18 vectors of norms 1 and 2 and partition them into three sets of size 6.
The resulting lattices have generator matrices
\beql{EqE6}
\left[ \begin{array}{ccr}
1 & 0 & 0 \\
0 & 1 & 1 \\
0 & 1 & -1
\end{array}
\right] , ~~
\left[ \begin{array}{ccr}
0 & 1 & 0 \\
1 & 0 & 1 \\
1 & 0 & -1
\end{array}\right] , ~~
\left[ \begin{array}{crc}
0 & 0 & 1 \\
1 & 1 & 0 \\
1 & -1 & 0
\end{array}
\right]
\eeq
and their intersection has generator matrix
\beql{EqE7}
\left[\begin{array}{ccc}
2 & 0 & 0 \\
0 & 2 & 0 \\
1 & 1 & 1
\end{array}
\right] ~,
\eeq
which is indeed another version of the bcc lattice.
This decomposition of the bcc lattice is the simplest we have found, and will be discussed further in Section 5.

For the fcc lattice generated by $(1,1,0)$, $(1,-1,0)$, $(0,1,-1)$, the vectors of small norms are:
$$
\begin{array}{ccccccc}
\multicolumn{3}{c}{\mbox{shape}} & ~~~~~ & \mbox{norm} & ~~~~~ & \mbox{number} \\ [+.05in]
0 & 0 & 0 && 0 && 1 \\
1 & 1 & 0 && 2 && 12 \\
2 & 0 & 0 && 4 && 6 \\
2 & 1 & 1 && 6 && 24 \\
2 & 2 & 0 && 8 && 12 \\
3 & 1 & 0 && 10 && 24 \\
2 & 2 & 2 && 12 && 8
\end{array}
$$
The simplest intersection we have found is formed by taking the 32
vectors
of norms 6 and 12 and partitioning them into four sets of size 8.
The resulting lattices have generator matrices
\beql{EqE8}
\left[\begin{array}{rrr}
2 & 1 & 1 \\
1 & 2 & -1 \\
-2 & 2 & 2
\end{array}
\right] , ~~
\left[ \begin{array}{rrr}
1 & 2 & 1 \\
-1 & 1 & 2 \\
2 & -2 & 2
\end{array}
\right] , ~~
\left[ \begin{array}{rrr}
1 & 1 & 2 \\
2 & -1 & 1 \\
2 & 2 & -2
\end{array}\right] , ~~
\left[ \begin{array}{rrr}
2 & -1 & -1 \\
-1 & 2 & -1 \\
2 & 2 & 2
\end{array}
\right]
\eeq
and their intersection is the fcc lattice with generator matrix
$(3,3,0)$, $(3,-3, 0)$, $(0,3,-3)$.
We will return to this example in Section 6.

\paragraph{Method 4: Using intersections of codes.}
A fourth method, which however has not yet led to any interesting
examples, is to reduce the problem to the analogous question for codes.
Let $\La (C)$ denote the lattice obtained by applying Construction A to a binary linear code $C$ (\cite[Chap. 5]{SPLAG}).
If $C_1, \ldots, C_r$ are codes of length $n$ whose intersection is a code $C$, then
$$\La (C) = \La (C_1) \cap \cdots \cap \La (C_r) ~.$$
This can be generalized to nonbinary codes \cite[Chaps. 7, 8]{SPLAG}, in particular to the case where the $C_i$ are nonbinary codes whose intersection is binary.

\section{The hexagonal lattice as an intersection of three lattices}
\hsp
We now begin our study of the honeycombs formed by some of the intersections described in Section 3.
Once the cells have been found it is generally straightforward to compute their volumes and second moments using the techniques presented in \cite{CSV} or \cite[Chap. 21]{SPLAG}.

The hexagonal lattice is the intersection of the three rectangular lattices given by \eqn{EqA1}, as in Fig.~\ref{A1}.
The honeycomb is shown in Fig.~\ref{A2}.
We now compute the mean squared error for this quantizer.
This analysis differs from that given
in \cite{BO1},
which used a different set of representation points for the quantizer.

There are four types of cells.
The origin is contained in a hexagon $\sP_1$ (horizontally shaded in Fig.~\ref{A2}) of edge length $1/\sqrt{3}$, area $V_1 = \sqrt{3}/2$ and second moment $U_1= 5 \sqrt{3} /72$.
The second type of cell is a small equilateral triangle $\sP_2$
(cross-hatched), with $V_2= \sqrt{3}/12$, $U_2 = \sqrt{3}/432$.
$\sP_3$ is an isosceles triangle (diagonally shaded), with $V_3 = \sqrt{3}/12$, $U_3 = 5 \sqrt{3} /1296$.
The fourth type, $\sP_4$, is a larger equilateral triangle (vertically shaded), with $V_4 = \sqrt{3}/4$, $U_4 = \sqrt{3}/48$.

\begin{figure}[htb]
\centerline{\psfig{file=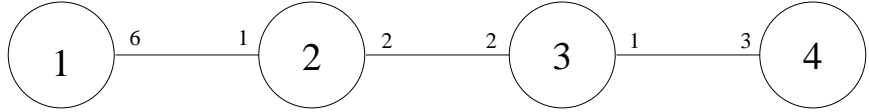,width=4.5in}}
\caption{Incidences between cells in honeycomb for hexagonal lattice.}
\label{IA2}
\end{figure}

The incidences between the different types of cells are shown in Fig.~\ref{IA2}.
Here (and in similar diagrams in later sections) a circle containing $i$ refers to a cell of type $\sP_i$, and an edge
\begin{center}
\psfig{file=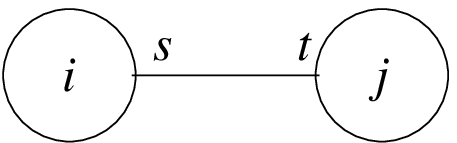,width=1.5in}
\end{center}
indicates that $\sP_i$ and $\sP_j$ share a common
face of the maximal dimension (here 1), and the edge labels indicate that each
$\sP_i$ is adjacent to $s$ cells of type $\sP_j$ and each $\sP_j$ to $t$ cells of type $\sP_i$.

The Voronoi cell $\sV$ for the intersection lattice is enclosed by the broken lines in Fig.~\ref{A2}.
In the notation of Section 2, $\sV$ contains $N_1=1$ copy of $\sP_1$, $N_2=6$ copies of $\sP_2$, $N_3=6$ copies of $\sP_3$ and $N_4=6 \times \frac{1}{3} =2$ copies of $\sP_4$, and the volume equation
\eqn{EqU0} reads
\beql{EqE9}
2 \sqrt{3} = \frac{\sqrt{3}}{2} + \frac{\sqrt{3}}{2} +
\frac{\sqrt{3}}{2} + \frac{\sqrt{3}}{2} ~.
\eeq
Thus the probabilities $p_1, \ldots, p_4$ of a randomly chosen point in the plane belonging to a cell of each type are all equal to 1/4.
From \eqn{EqU5}, the mean squared error for this quantizer is $G= 2^{7/4} /27 = 0.1246 \ldots$
This value (and corresponding values of $G$ found
in the next three Sections) is considerably worse than the value $0.080188 \ldots$ for the hexagonal lattice itself.

\section{The bcc lattice as an intersection of three lattices}
\hsp
The bcc lattice is the intersection of the three ``rectangular'' lattices
$\La_1, \La_2 , \La_3$ defined in \eqn{EqE6}.
Each $\La_i$ is congruent to $\ZZ \times \sqrt{2} \ZZ \times \sqrt{2} \ZZ$, and has as Voronoi cell a brick with square cross-section.
There is an obvious symmetric group $S_3$ that permutes the $\La_i$.

Again the honeycomb contains four types of cells.
The origin is contained in the intersection of the Voronoi cells at 0 for the three $\La_i$.
This is a cube, $\sP_1$, with vertices $(\pm 1/2, \pm 1/2, \pm 1/2)$, volume $V_1 =1$, and second moment $U_1 = 1/4$.
Across each square face of $\sP_1$ is a square pyramid $\sP_2$, such as that with base $(1/2, \pm 1/2 , \pm 1/2 )$, apex $(1,0,0)$,
$V_2 =1/6$, $U_2 = 11/600$.
Across each triangular face of $\sP_2$ is a tetrahedron $\sP_3$, such as that with vertices
$(1/2, \pm 1/2 , 1/2)$, $(1,0,0)$, $(1,0, 1/2)$, $V_3 = 1/24$, $U_3=1/512$.
Finally, across the other three faces of $\sP_3$ we reach a fourth type of cell, $\sP_4$, a quarter-octahedron, which occurs in two orientations, one having vertices such as
\beql{Eq5A}
\left( \frac{1}{2}, 0, \frac{1}{2} \right),~~
\left( \frac{1}{2}, 0,1 \right), ~~
\left( 1,0, \frac{1}{2} \right) , ~~(1,0,1) , ~~
\left( \frac{1}{2} , \pm \frac{1}{2} , \frac{1}{2} \right),
\eeq
the other having vertices such as
\beql{Eq5B}
\left( 1, \frac{1}{2}, \frac{1}{2} \right) ,~~
\left(1, 0, \frac{1}{2} \right),~~
\left( 1, \frac{1}{2}, 0 \right), ~~(1,0,0), \left(
\frac{1}{2}, \frac{1}{2}, \frac{1}{2} \right) ,
\left( \frac{3}{2} , \frac{1}{2} \frac{1}{2} \right) \,.
\eeq
These may be described as quarters of squat octahedra.
E.g., \eqn{Eq5A} is a quarter of the octahedron with vertices
$(\pm 1, 0, \pm 1)$, $(1/2 , \pm 1/2, 1/2)$.
For $\sP_4$ we have $V_4 = 1/12$, $U_4 = 1/192$.

\begin{figure}[htb]
\centerline{\psfig{file=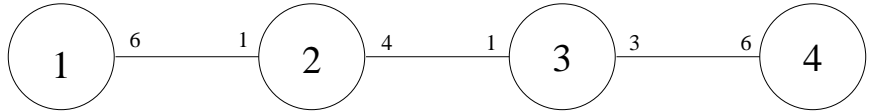,width=4.5in}}
\caption{Incidences between cells for honeycomb for bcc lattice.}
\label{IBCC}
\end{figure}

No further types of cell appear:
every face of either version of $\sP_4$ leads to a $\sP_3$.
The incidence diagram is shown in Fig.~\ref{IBCC}.

The Voronoi cell $\sV$ for the
intersection lattice is a truncated octahedron with 24 vertices $(0, \pm 1/2 , \pm 1)$.
This contains $N_1=1$ copy of $\sP_1$, $N_2=6$ copies of $\sP_2$ and $N_3=24$ copies of $\sP_3$.
The cells of type $\sP_4$ partially overlap $\sV$.
There are 12 of type \eqn{Eq5A}, intersecting $\sV$ in a tetrahedron such as that with vertices $(1/2, \pm 1/2 , 1/2 )$, $(1/2, 0,1)$, $(1,0,1/2)$, with volume $1/24$.
There are also 24 of type \eqn{Eq5B}, intersecting $\sV$ in a tetrahedron such as
$(1/2, 1/2, 1/2)$, $(1,0,0)$, $(1,0, 1/2)$, $(1,1/2, 0)$, with volume $1/48$.
The volume equation \eqn{EqU0} reads
\begin{eqnarray}\label{Eq5C}
4 & = & 1 \times 1 + 6 \times \frac{1}{6} + 24 \times \frac{1}{24} + 12 \times \frac{1}{24} + 24 \times \frac{1}{48} \nonumber \\
& = & 1 + 1 + 1 + 1 ~,
\end{eqnarray}
so again the probabilities $p_i$ of a randomly chosen point belonging to a cell of given type are all equal to $1/4$.
The mean squared error is
$$G = \frac{751 \sqrt{3}}{9600} = 0.1355 \ldots$$

\begin{figure}[htb]
\centerline{\psfig{file=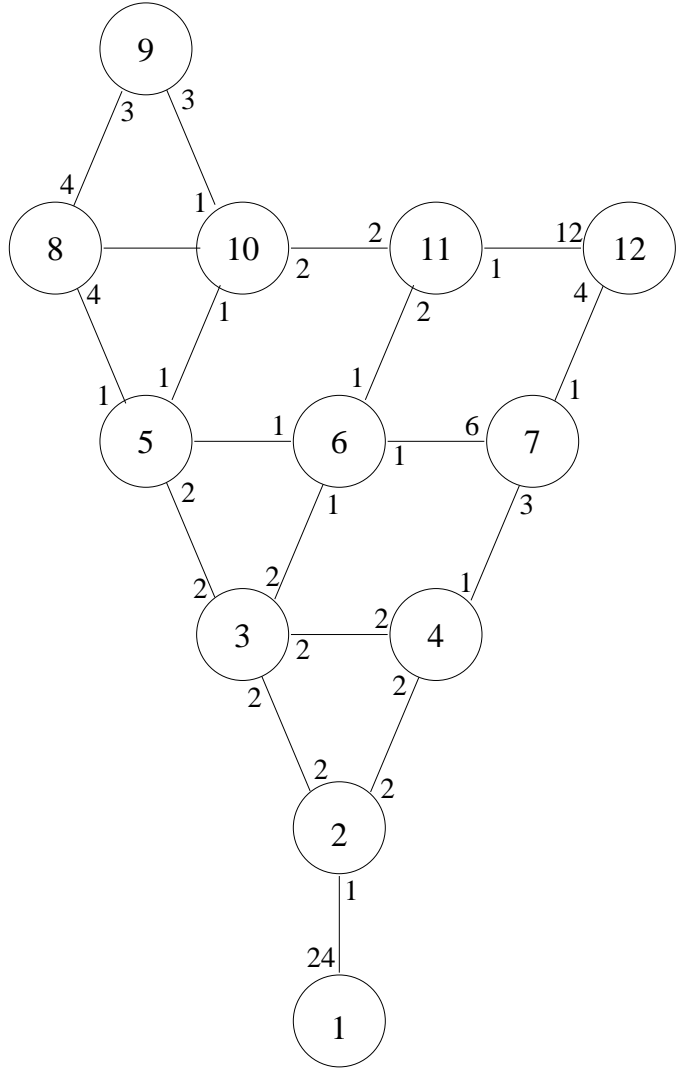,width=3.5in}}
\caption{Incidences among cells of fcc honeycomb.}
\label{ID3}
\end{figure}

\section{The fcc lattice as an intersection of four lattices}
\hsp
The fcc lattice is the intersection of the four lattices $\La_1, \ldots, \La_4$ defined in \eqn{EqE8}.
Each of these has Gram matrix equivalent to
$$
\left[
\begin{array}{rrr}
6 & -3 & 0 \\
-3 & 6 & 0 \\
0 & 0 & 12
\end{array}
\right]
$$
and is a direct sum $\sqrt{3} A_2 \oplus \sqrt{12} \ZZ$, with Voronoi cell a hexagonal prism.
The $\La_i$ look more symmetrical if they are written in the coordinates used to describe
the root lattice $A_3$
$(\cong D_3 )$, that is, using four coordinates that add to 0.
Then $\La_1$ has generator matrix
$$
\left[
\begin{array}{crrr}
0 & 2 & -1 & -1 \\
0 & -1 & 2 & -1 \\
3 & -1 & -1 & -1
\end{array}
\right]
$$
and the others are given by cyclic shifts of these columns.
This shows that there is a symmetric group $S_4$ permuting the $\La_i$.
However, the three-dimensional coordinates given in \eqn{EqE8} are more convenient for computations.

\begin{figure}[htb]
\centerline{\psfig{file=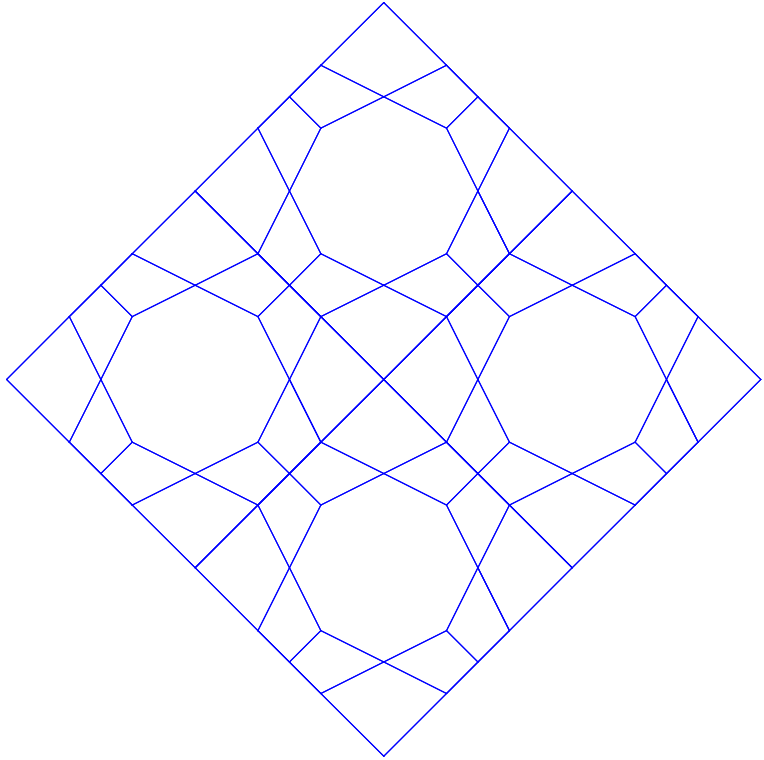,width=4in}}
\caption{Cross-section of fcc honeycomb along plane $z=0$ $(- 3 \le x \le 9$, $-6 \le y \le 6)$, with origin at center of octagon on left.
Only three cells are visible: $\sP_1$ (octagon), $\sP_2$ (small kite), $\sP_9$ (large dart).}
\label{HD1}
\end{figure}

\begin{figure}[htb]
\centerline{\psfig{file=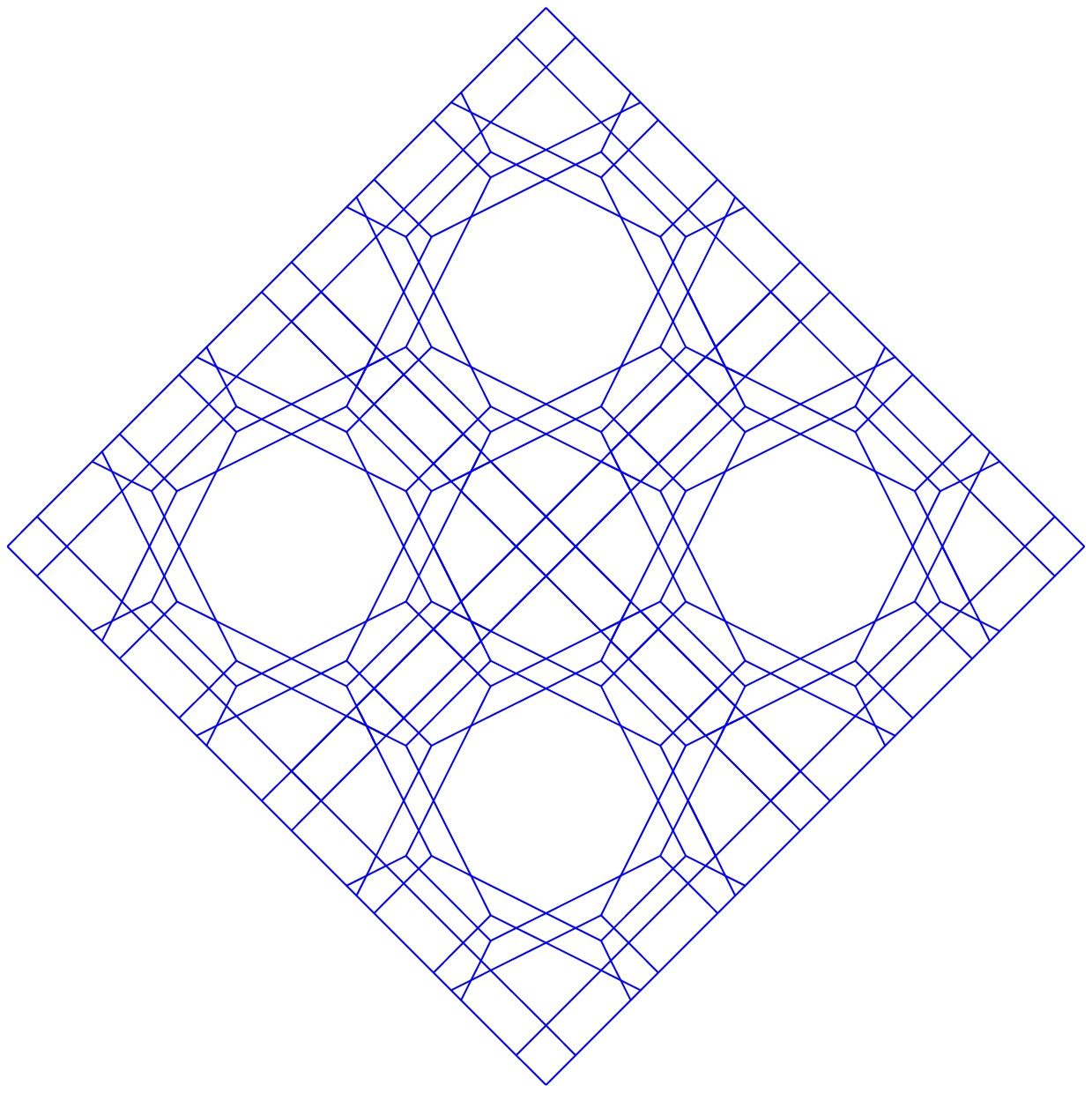,width=4in}}
\caption{Cross-section through fcc honeycomb along plane $z= 0.35$.  Cross-sections
of all 12 types of cells can be seen.}
\label{HD2}
\end{figure}

This honeycomb is the most complicated we have analyzed and we shall give only a brief description.
There are twelve types of cells,
$\sP_1, \ldots, \sP_{12}$, whose parameters are summarized in Table~\ref{TD31} and whose incidences are shown in Fig. \ref{ID3}.
Figures \ref{HD1} and \ref{HD2} shows cross-sections through the honeycomb along the planes $z=0$ and $z=0.35$.
\begin{table}[htb]
$$
\begin{array}{cccccccc}
i & v & e & f & N_i & V_i & p_i & U_i \\ [+.05in]
1 & 26 & 48 & 24 & 1 & 9 & 1/6 & 1449/160 \\
2 & 5 & 8 & 5 & 24 & 1/8 & 1/18 & 41/3200 \\
3 & 6 & 12 & 8 & 24 & 1/4 & 1/9 & 9/320 \\
4 & 5 & 8 & 5 & 24 & 1/8 & 1/18 & 41/3200 \\
5 & 5 & 9 & 6 & 24 & 1/10 & 2/45 & 427/50000 \\
6 & 4 & 6 & 4 & 48 & 1/40 & 1/45 & 443/320000 \\
7 & 10 & 18 & 10 & 8 & 27/40 & 1/10 & 41013/160000 \\
8 & 6 & 12 & 8 & 6 & 1 & 1/9 & 47/180 \\
9 & 5 & 9 & 6 & 8 & 3/8 & 1/18 & 189/3200 \\
10 & 4 & 6 & 4 & 24 & 1/40 & 1/90 & 1897/1280000 \\
11 & 5 & 8 & 5 & 24 & 3/40 & 1/30 & 2319/400000 \\
12 & 22 & 36 & 16 & 2 & 63/10 & 7/30 & 213597/40000
\end{array}
$$
\caption{The twelve types of cells in the fcc honeycomb, showing numbers of vertices, edges, faces $(v,e,f)$, the number per fcc cell $(N_i)$, and their volumes, probabilities and second moments$(V_i, p_i, U_i )$.}
\label{TD31}
\end{table}

The following is a brief description of the cells, including coordinates for one cell of each type.

$\sP_1$.
Obtained from cube by pushing in corners and pulling out edge midpoints.
Vertices: all cyclic shifts and sign changes of $(3/2, 0,0)$, $(1,1,0)$, $(3/4, 3/4, 3/4)$.

$\sP_2$.
Pyramid with kite-shaped base.
Base: $(3/2, 0,0)$, $(1,1,0)$, $(1,0,1)$, $(3/4, 3/4, 3/4)$,
apex: $(3/2, 1/2, 1/2)$.

$\sP_3$.
Irregular octahedron with vertices $(3/2, 0,0)$, $(3/2, 3/2, 0)$,
$(1,1,0)$, $(2,1,0)$, $(3/2, 1/2, \pm 1/2)$.

$\sP_4$.
Congruent to $\sP_2$.
Vertices:
$(3/2, 0,0)$,
$(3/4, 3/4, 3/4)$, $(3/2, 1/2, 1/2)$, $(1/2, 3/2, 1/2)$, $(1,1,0)$.

$\sP_5$.
Irregular polyhedron with six faces.
Vertices:
$(3/2, 0,0)$, $(2,1,0)$, $(2,0,1)$,
$(9/5, 3/5, 3/5)$, $(3/2, 1/2, 1/2)$.

$\sP_6$.
Tetrahedron:
$(2,1,0)$, $(3/2, 3/2, 0)$, $(3/2, 1/2, 1/2)$,
$(9/5, 3/5, 3/5)$.

$\sP_7$.
``Flying saucer'':
hexagonal base (cyclic shifts of $(3/2, 3/2,0)$ and $(9/5, 3/5, 3/5)$) with four vertices above it (cyclic shifts of $(3/2, 1/2, 1/2)$ and
$(3/4, 3/4, 3/4)$ at apex).

$\sP_8$.
Irregular octahedron with vertices $(2, \pm 1, 0)$, $(2, 0, \pm 1)$,
$(3/2, 0,0)$, $(3,0,0)$.

$\sP_9$.
Regular tetrahedron (vertices $(2,1,0)$,
$(2, 0,1)$, $(3,0,0)$,
$(3,1,1)$)
with triangular cap (apex $(9/4, 3/4, 3/4)$) on one face.

$\sP_{10}$.
Irregular tetrahedron with vertices $(2,1,0)$, $(2,0,1)$,
$(9/4, 3/4, 3/4)$, $(9/5, 3/5, 3/5)$.

$\sP_{11}$.
Another pyramid with kite-shaped base.
Base:
$(3/2, 3/2, 0)$, $(9/4, 3/4, 3/4 )$, $(9/5, 3/5, 3/5)$, $(12/5, 6/5, 3/5)$,
apex: $(2,1,0)$.

\begin{figure}[htb]
\centerline{\psfig{file=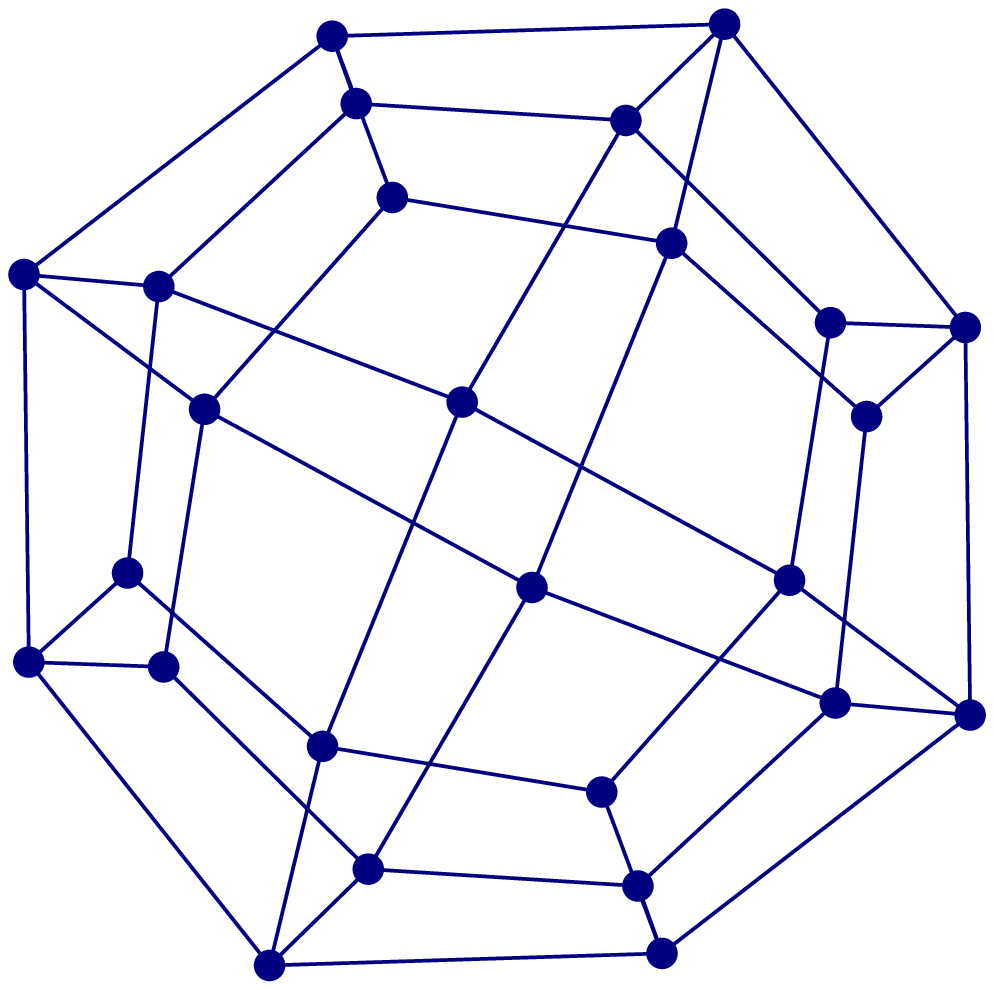,width=4in}}
\caption{26-vertex cell $\sP_{12}$.}
\label{HD3}
\end{figure}

$\sP_{12}$.
See Fig. \ref{HD3}.
Has 26 vertices, four hexagonal and 12 kite-shaped faces.
Vertices are all permutations of
$(3/2, 3/2, 0)$, $(3/2, 3/2, 3)$, $(9/4, 3/4, 3/4)$,
$(9/4, 9/4, 9/4)$, $(9/5, 3/5, 3/5)$,
$(12/5, 9/5, 9/5)$, $(12/5, 6/5, 3/5)$.

If we decompose $\RR^3$ into Voronoi cells for the intersection lattice $\La$ \eqn{EqE9}, just three of the twelve types of cells are cut by the boundary walls.
Cells of type $\sP_9$ are cut into three equal pieces, cells of type $\sP_{11}$ are cut in half, and cells of type $\sP_{12}$ are cut into four equal ice-cream cone shaped pieces.
The base of each cone is at the center of $\sP_{12}$, and the top contains one of the hexagons and parts of the neighboring faces.

The mean squared error is
$$G ~=~ \frac{12269777}{816480000}~ 3^{28/135}~ 5^{8/27}~ 7^{38/45} ~=~ 0.1572 \dots ~.$$

\section{The $D_4$ lattice as an intersection of three cubic lattices}
\begin{figure}[htb]
\centerline{\psfig{file=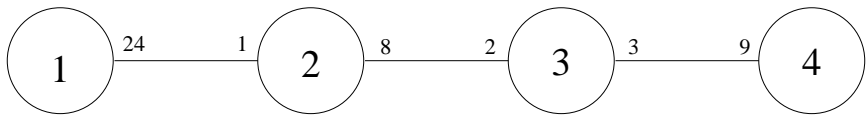,width=4.5in}}
\caption{Incidences among cells of $D_4$ honeycomb.}
\label{ID4}
\end{figure}

\hsp
This is the nicest example we have found.
The three cubic lattices
$\La_1$, $\La_2$, $\La_3$ are defined in \eqn{EqD1} and their intersection $\La \cong D_4$ in \eqn{EqD2}.
There are just four types of cells, 
$\sP_1$, $\sP_2$, $\sP_3$, $\sP_4$, whose properties are summarized in Table \ref{TD41} and whose intersections are shown in Fig.~\ref{ID4}.
We use coordinates $(a,b,c,d)$ for points in $\RR^4$.
\begin{table}[htb]
$$\begin{array}{ccccccc}
i & v & f & N_i & V_i & p_i & U_i \\ [+.05in]
1 & 24 & 24 & 1 & 8 & 1/4 & 104/15 \\
2 & 7 & 9 & 24 & 1/3 & 1/4 & 8/105 \\
3 & 5 & 5 & 96 & 1/12 & 1/4 & 11/900 \\
4 & 6 & 9 & 32 & 1/4 & 1/4 & 1/20
\end{array}
$$
\caption{Cells in $D_4$ honeycomb, showing numbers of vertices and 3-dimensional faces $(v,f)$, the number per $D_4$ cell $(N_i)$ and their volumes, probabilities and second moments $(V_i, p_i, U_i )$.}
\label{TD41}
\end{table}

$\sP_1$ is the 4-dimensional regular polytope known as a 24-cell
(\cite{SPLAG}, \cite{Cox}).
The Voronoi cells for $\La_1$, $\La_2$, $\La_3$ at the origin are all cubes, whose intersection is bounded by the hyperplanes
$$|a| \le 1,~
|b| \le 1,~
|c| \le 1, ~
|d| \le 1, ~
|a|+|b| + |c| +|d| \le 1 ,
$$
which is the 24-cell with vertices of the form $(\pm 1, \pm 1, 0,0)$.

Across each of the 24 octahedral faces of $\sP_1$ we reach an octahedral-based pyramid $\sP_2$, such as that with vertices $(1, \pm 1, 0,0)$, $(1,0, \pm 1, 0)$, $(1,0,0,\pm 1)$ and $(2,0,0,0)$ (the apex).

There are eight other faces of $\sP_2$, regular tetrahedra;
these lead to copies of $\sP_3$,
which is an irregular simplex such as that with vertices
$(1,1,0,0)$, $(1,0,1,0)$, $(1,0,0,1)$,
$(2,0,0,0)$, $(1,1,1,1)$.

Finally, three of the five faces of each $\sP_3$ lead to cells of the fourth type, $\sP_4$.
This can best be described as the product of two skew equilateral triangles of different sizes (just as a tetrahedron in three dimensions is the product of two skew line segments).
Take an equilateral triangle with vertices $p= (2,1,1,0)$, $q= (1,1,0,0)$, $r= (1,0,1,0)$ and another with vertices $P= (2,0,0,0)$,
$Q= (1,1,1,1)$, $R= (1,1,1,-1)$.
Then $\sP_4$ is their convex hull.
There are nine tetrahedral faces given by the convex hull of an edge of the first triangle and an edge of the second triangle.

If we decompose $\RR^4$ into Voronoi cells for the intersection lattice, only cells of type $\sP_4$ are cut by the boundary walls.
Each $\sP_4$ is divided into three equal pieces, a typical piece being an
``ice-cream cone'' whose center is at the center of $\sP_4$ and whose three-dimensional face is the convex hull of the second triangle $(P,Q,R)$ and any edge of the first triangle.

The volume equation \eqn{EqU0} then reads
$$
32 = 1 \times 8 + 24 \times \frac{1}{3} + 96
\times \frac{1}{12} + 96 \times \frac{1}{4} \times \frac{1}{3} \,.
$$
Again a random point is equally likely to fall into a cell of any of the four types.
The mean squared error is
$$G ~=~ \frac{757}{8400}~ 2^{1/8}~ 3^{1/4} ~=~ 0.1293 \dots ~.$$

\section{Conclusions and comments}
\hsp
In each of the honeycombs of Sections 4, 5 and 7 just four types of cells occurred.
This is easily explained in the case of the $D_4$ honeycomb:
there are three equivalent lattices $\La_1$, $\La_2$, $\La_3$, and the
associated quantizers are essentially making binary decisions about the location of a point with respect to the intersection lattice.
So the space is divided up into regions that can be labeled 000, 001, 011 and 111.

This argument doesn't quite apply to the $A_2$ or bcc honeycombs, since there the individual lattices themselves are not fully symmetric (rectangular rather square in the $A_2$ case, for instance).
So it is fortuitous that only four cells occur.

In contrast, the fcc honeycomb of Section 6 shows that the number of cells can increase rapidly in less fortunate cases with more component lattices.
It would be interesting to see how complicated are the $E_6^\ast$ and $E_8$ honeycombs mentioned in Section 3.

As for practical applications, it is possible that better quantizers could be obtained by amalgamating less symmetrical cells.
For example, in Fig. \ref{A2}, the diagonally and vertically shaded triangles 
could be amalgamated to give a honeycomb made up of regular hexagons and equilateral triangles with the same edge length as the hexagons.
The new honeycomb will have larger absolute error but a smaller normalized error $G$.
We did not investigate this possibility for the higher-dimensional examples.

Another topic for future research is the design of quantizers by taking a simple
initial lattice $\La_1$ and combining it with several other lattices
which are both rotations {\em and translations} of $\La_1$.

\subsection*{Acknowledgements}
\hsp
We thank John Conway, Suhas Diggavi, Eric Rains and Vinay Vaishampayan for helpful discussions.
Deborah Swayne provided valuable help in running XGobi, and David Johnson kindly supplied the graph coloring program that was used in Section 3.

\end{document}